\documentclass[10pt]{amsart}

\usepackage{amsfonts}
\usepackage{amssymb}
\usepackage{graphics}
\usepackage{amsmath, amsthm, latexsym}

\usepackage{hyperref}

\usepackage[margin=1.2in]{geometry}

\def \c{\mathbb{C}}
\def \z{\mathbb{Z}}
\def \r{\mathbb{R}}
\def \n{\mathbb{N}}
\def \p{\mathbb{P}}
\def \q{\mathbb{Q}}

\def \W{\mathcal{W}}

\def \t{\mathfrak{t}}
\def \k{{\bf k}}

\def \Q{\mathcal{Q}}

\def \SL{\textup{SL}}

\def \GL{\textup{GL}}
\def \PGL{\textup{PGL}}
\def \PO{\textup{PO}}
\def \SP{\textup{Sp}}
\def \SO{\textup{SO}}
\def \Pic{\textup{Pic}}

\def \B{\mathcal{B}}
\def \s{{\underline{w}_0}}
\def \.{\cdot}

\def \deg{\textup{deg}}
\def \dim{\textup{dim}}

\def \Vol{\textup{Vol}}

\def \Sym{\textup{Sym}}
\def \Lie{\textup{Lie}}
\def \ker{\textup{ker}}

\def \Pic{\textup{Pic}}
\def \End{\textup{End}}

\theoremstyle{plain}
\newtheorem{Th}{Theorem}[section]
\newtheorem{Lem}[Th]{Lemma}

\newtheorem{Cor}[Th]{Corollary}

\theoremstyle{definition}
\newtheorem{Ex}[Th]{Example}
\newtheorem{Def}[Th]{Definition}
\newtheorem{Rem}[Th]{Remark}

\begin{document}
\title{Note on cohomology rings of spherical varieties and volume polynomial}
\author{Kiumars Kaveh}

\begin{abstract}
Let $G$ be a complex reductive group and $X$ a projective
spherical $G$-variety. Moreover, assume that the subalgebra $A$ of the
cohomology ring $H^*(X, \r)$ generated by the Chern classes of line bundles has Poincare
duality. We give a description of the subalgebra $A$ in terms of
the volume of polytopes. This generalizes the Khovanskii-Pukhlikov
description of the cohomology ring of a smooth toric variety.
In particular, we obtain a unified
description for the cohomology rings of complete flag varieties and smooth toric varieties.
{As another example we get a description of the cohomology ring of the variety of complete conics.}
We also address the
question of additivity of the moment and string polytopes and prove
the additivity of the moment polytope for {complete} symmetric varieties.
\end{abstract}

\maketitle


%\tableofcontents

\section*{Introduction} \label{sec-intro}
The well-known theorem of Kushnirenko
(\cite{Bernstein}, \cite{Kushnirenko}, \cite{BKK})
on the number of solutions of a system of Laurent polynomial equations
in terms of the volumes of their Newton polytopes, can be
interpreted as giving a formula for the intersection numbers of
divisors in a projective toric variety. Khovanskii and Pukhlikov
\cite[$\S$1.4]{Kh-P} observed that the Kushnirenko
theorem, in fact, completely determines the cohomology ring of a
smooth projective toric variety. The key fact is that the cohomology
ring of a projective toric variety is generated in degree
$2$.

In this note we apply the Khovanskii-Pukhlikov approach to the cohomology ring of
projective spherical varieties.

The Kushnirenko theorem has been generalized to
representations of reductive groups by B. Kazarnovskii
\cite{Kazarnovskii} and more generally, to spherical varieties by M.
Brion \cite[$\S$ 4.1]{Brion2}. Let $G$ be a complex connected
reductive algebraic group. A normal $G$-variety $X$ is called {\it spherical},
if a Borel subgroup of $G$ has a dense open orbit. Equivalently, $X$ is
spherical if the space of global sections of any $G$-linearized line bundle
on $X$, regarded as a $G$-module, is multiplicity free. Let $V$ be a
finite dimensional $G$-module and $X \hookrightarrow \p(V)$ a closed
spherical $G$-subvariety of dimension $n$. Let $L$ be the line bundle on
$X$ obtained by the restriction of anti-canonical line bundle
$\mathcal{O}_{\p(V)}(1)$ to $X$. It is naturally a $G$-linearized line bundle. 
The Brion-Kazarnovskii formula gives the
degree of $X$ in $\p(V)$ as the integral of an explicitly defined
polynomial over the {\it moment polytope} $\mu(X, L)$ (Theorem
\ref{thm-Brion}). In terms of $L$, the degree
of $X$ in $\p(V)$ is the self-intersection number of the Chern class $c_1(L)$.

In order to answer a question posed by A. G. Khovanskii, when $G$ is one of the classical
groups over $\c$, A. Okounkov \cite{Okounkov} associated a polytope $\Delta(X, L) \subset \r^n$
to the pair $(X, L)$ such that $\deg(X)$ in $\p(V)$ is equal to $n!\Vol_n(\Delta(X, L))$ (where $n
= \dim(X)$ and $\Vol_n$ is the usual Lebesgue measure in $\r^n$).
The polytope $\Delta(X, L)$ is the polytope fibred over the
moment polytope of $(X, L)$ with the Gelfand-Cetlin polytopes as fibers.
Later, following ideas of Caldero \cite{Caldero}, Alexeev and Brion observed that
Okounkov's construction can be done for spherical varieties of any
connected reductive group \cite{Alexeev-Brion}. They used it to extend the
toric degeneration result in \cite{Kaveh} for $G=\SP(2n, \c)$, to any
connected reductive group $G$. That is, $X$ can be deformed, in a flat family, to the toric
variety associated to any of these polytopes.
Alexeev-Brion-Caldero's
work relies on the generalization of the well-known Gelfand-Cetlin
polytopes for $\GL(n, \c)$ to any complex connected reductive group. These are the so-called
{\it string polytopes}. Recall that a
reduced decomposition $\s$ for $w_0$, the longest element in the
Weyl group, is a decomposition $w_0 = s_{i_1}\cdots s_{i_N}$ into simple
reflections. Here $N = \ell(w_0)$ is the length of $w_0$ which is
equal to the number of positive roots. For a choice of a reduced
decomposition for $w_0$ and a dominant weight $\lambda$, one
constructs a string polytope $\Delta_{\s}(\lambda) \subset \r^N$
\cite{Littelmann}. The string polytope has the property that the
integral points inside it are in one-to-one correspondence with the
elements of a so-called crystal basis for the highest weight
representation space $V_\lambda$. For a spherical variety $X$ and a
$G$-linearized very ample line bundle $L$ on $X$, the polytope constructed by Alexeev-Brion-Okounkov
is a rational polytope $\Delta_{\s}(X, L) \subset \Lambda_\r \times \r^N$, where
$\Lambda_\r$ is the real vector space generated by the weight
lattice. The projection on the first factor maps this polytope onto
the moment polytope $\mu(X, L)$ and the fibre over each $\lambda \in
\mu(X,L)$ is the string polytope $\Delta_{\s}(\lambda)$. From the
above mentioned property of the string polytopes, and using the asymptotic Riemann-Roch,
it can be seen that
\begin{equation} \label{eqn-deg-introduction}
\deg(X, L) = n! \Vol_n(\Delta_{\s}(X, L)).
\end{equation}
In analogy with toric geometry, we call the polytope
$\Delta_{\s}(X, L)$, the {\it Newton polytope of $(X, L)$}.

{It follows from the theory of Bialynicki-Birula that a smooth spherical variety $X$
has a paving by (complex) affine spaces}. Hence there is a ring isomorphism between 
$H^*(X, \z)$ and the Chow ring $CH^*(X)$. Then {the Picard group $\Pic(X)$ can be identified with 
$H^2(X, \z)$ and the real span of $\Pic(X)$ is $H^2(X, \r)$}. 

Unlike the toric case, the cohomology ring of a
projective spherical variety, in general, {is not generated in degree $2$, and in particular
by the Chern classes of line bundles} (for example the Grassmannian of $2$-planes in $\c^4$).
Let $A$ denote the subalgebra of the cohomology ring of a
projective spherical variety $X$ generated by {the Chern classes of line bundles}.
In this note, we use the Alexeev-Brion-Okounkov formula  (\ref{eqn-deg-introduction})
for the degree, to give a description of the graded algebra $A$,
provided that this subalgebra satisfies Poincare duality.
This description is in terms of volumes of Newton polytopes (Theorem \ref{thm-main}). It is a direct
generalization of the Khovanskii-Pukhlikov description of the cohomology
ring of a smooth projective toric variety (\cite{Kh-P}). It will follow from a general theorem
in commutative algebra which asserts that certain algebras can be
realized as a quotient of the algebra of differential operators with constant coefficients
(Theorem \ref{thm-comm-alg}). {As far as the author knows there are no examples 
of a variety known such that the subalgebra generated by the Chern classes of line bundles does not satisfy 
Poincare duality}.

It is well-known that the cohomology ring of the complete flag
variety $G/B$ is generated by {the Chern classes of line bundles}. Thus our result in particular gives a
description of the cohomology ring of $G/B$ in terms of
the volumes of polytopes (Corollary \ref{cor-G/B}). It is
interesting to note that this way, we arrive at analogous
descriptions for the cohomology rings of toric varieties and the
flag variety $G/B$. This description of the cohomology ring of $G/B$
coincides with Borel's description via a theorem of Kostant
\cite{Kostant}. It states that if $P$ is the polynomial which is the
product of equations of the hyperplanes perpendicular to the roots,
then the ideal of differential operators annihilating $P$ is
generated by the $W$-invariant operators.

{As another example 
we get a description of the cohomology ring of the variety of complete conics in terms of volume 
of polytopes. It is not hard to see that this description coincides with the classical description of this cohomology 
in terms of generators and relations (Section \ref{ex-complete-conics}).}

Let $X$ be a projective $T$-toric variety and $L$ a very ample $T$-linearized line bundle on $X$.
A nice property of the Newton polytope map $L \mapsto \Delta(X, L)$ in the toric case
is the additivity. Namely if $L_1, L_2$ are two very ample $T$-linearized line bundles
then $$\Delta(X, L_1 \otimes L_2) = \Delta(X, L_1) + \Delta(X, L_2),$$
where the addition is the Minkowski sum of convex polytopes. 
This implies that the formula $\deg(X, L) = n! \Vol_n(\Delta(X, L))$ extends to give a
formula for the product of Chern classes of $n$ line bundles as the mixed volume of
the corresponding Newton polytopes. This is the content of the Bernstein theorem in toric geometry.

Now let $X$ be a projective spherical $G$-variety.
To investigate the analogy with the toric case more closely, in
Section \ref{sec-additive} we address the question when the map
$$L \mapsto \Delta_\s(X, L)$$ is additive i.e. $$\Delta_\s(X, L_1
\otimes L_2) = \Delta_\s(X, L_1) + \Delta_\s(X, L_2),$$ where $L_1$
and $L_2$ are very ample $G$-linearized line bundles. For this, we discuss the
additivity of the moment and string polytopes separately. We will see
that in general neither one is additive. Although we show that when
$X$ is a {\it complete symmetric variety} and the restriction of $L$ to
the open orbit is trivial, then the moment polytope map is
additive (Theorem \ref{th-symm-var-additive}), that is $$\mu(X, L_1 \otimes L_2) = \mu(X, L_1) + \mu(X, L_2).$$

{Also there is an important special class of string polytopes, namely the Gelfand-Cetlin polytopes, 
which are additive: For a (complex) connected reductive 
group $G$ whose Lie algebra is a direct sum of classical simple Lie algebras and/or a commutative Lie algebra, and $\lambda \in \Lambda^+_\r$, one defines the Gelfand-Cetlin, or for short G-C, polytope
$\Delta(\lambda)$ (see \cite{B-Z}). One can write down explicitly the defining equations of $\Delta(\lambda)$
and from this it easily follows that $\lambda \mapsto \Delta(\lambda)$ is additive, i.e., 
for $\lambda_1$, $\lambda_2 \in \Lambda^+_\r$:
$$\Delta(\lambda_1 + \lambda_2) = \Delta(\lambda_1) + \Delta(\lambda_2).$$
From this the additivity of $\Delta_\s(X,L)$ follows for complete symmetric varieties of such a group $G$ and 
a reduced decomposition $\s$ which gives rise to Gelfand-Cetlin polytopes. 
An interesting class of examples of symmetric varieties are compactifications of reductive groups.}

{As in the toric case, the additivity of the Newton polytope then implies a 
Bernstein theorem for complete symmetric varieties. That is, the 
intersection number of the Chern classes of $G$-linearized 
line bundles can be obtained 
from the mixed volume of the corresponding Newton polytopes (Theorem \ref{th-Bernstein-for-symm-var}).
The answer to several well-known 
classical problems in enumerative geometry can be represented as the intersection 
numbers of divisors in certain symmetric varieties (\cite[Section 10]{DP}). 
The above Bernstein theorem for symmetric varieties can be applied to represent the answers 
as mixed volumes of Newton polytopes, which then can be computed explicitly. As an example, 
we expect that in this fashion one can recover the answer to a classical problem of Schubert (Remark \ref{rem-enumerative-application}).}

{In general, for a projective spherical variety $X$, we can show that there is a subdivision of the cone of ample $G$-linearized line bundles into smaller cones such that the map $L \mapsto \mu(X, L)$ is additive on each cone of this subdivision.
The same is true for the string polytopes, i.e. there is a subdivision of the positive Weyl chamber into smaller cones such that $\lambda \mapsto \Delta_\s(\lambda)$ is additive on each cone.}

In \cite{Caldero}, Caldero gives a flat degeneration of the flag variety $G/B$ together with a 
line bundle $L_\lambda$ (corresponding to a dominant weight $\lambda$)
to the toric variety $X_0$ {together with a $\q$-divisor class $L_0$ corresponding to the polytope
$\Delta_{\s}(\lambda)$. (In general, the polytope $\Delta_\s(\lambda)$ may not be integral and hence 
determines a toric variety together with a $\q$-divisor class.)} Since the Hilbert polynomial is preserved in a
flat family, {the degree of $L_\lambda$ and the self-intersection of the $\q$-divisor $L_0$ are
the same}. This provides a geometric visualization for the similarity of our description of
the cohomology rings of $G/B$ and toric varieties. Note that in general the cohomology ring is not
preserved under a flat degeneration, {in fact the cohomology ring of the degeneration $X_0$ is usually
much larger than that of $X$}. Also the toric variety corresponding to a string polytope
(in particular a G-C polytope) typically is not smooth.
As mentioned above, for a spherical variety $X$ and a
very ample $G$-linearized line bundle $L$, the toric degeneration to the toric variety
associated to the polytope $\Delta_{\s}(X, L)$ has been constructed in \cite{Kaveh} (in the special case of $G=\SP(2n, \c)$) and \cite{Alexeev-Brion} (in general).

Recently the approach of the present paper, in particular the description of the
cohomology ring of the flag variety in terms of volumes of the G-C polytopes, has been used
in \cite{Valentina}. In this work the author
develops a remarkable connection between the Schubert calculus (in type $A$) and
the combinatorics of the G-C polytopes.

The present note is also closely related to \cite{G-S}.

Few words about the organization of the paper: Section \ref{sec-comm-alg}
discusses a commutative algebra theorem which realizes certain algebras as
quotients of the algebra of differential operators with
constant coefficients.
%Section \ref{sec-Poincare-duality} shows that the
%subalgebra of cohomology generated by $\Pic(X)$ satisfies the conditions of this theorem.
Section \ref{sec-degree} discusses the Brion-Kazarnovskii formula for
the degree of a projective spherical variety,
string polytopes and Newton polytopes for spherical varieties. Section \ref{sec-additive} addresses the
question of additivity of the moment and string polytopes and proves the additivity of the
moment polytope for complete symmetric varieties. The main theorem appears in Section \ref{sec-main}
and finally Section \ref{sec-examples} discusses the examples of toric varieties, flag varieties
{and the variety of complete conics}.\\

\noindent{\bf Acknowledgment:} The author would like to thank M. Brion,
A. G. Khovanskii, V. Kiritchenko, E. Meinrenken and V. Timorin for
helpful discussions, comments and suggestions. In particular M. Brion suggested a 
counter-example in Section \ref{sec-additive} and also suggested to extend
the additivity of moment polytope to symmetric varieties.
Also thanks to M. Brion and W. Fulton for pointing out some errors in the previous versions.
The application of Theorem \ref{thm-main} to the variety of complete conics was suggested to the author
by V. Kiritchenko.

\section{A theorem from commutative algebra} \label{sec-comm-alg}
Let $A$ be a finite dimensional graded commutative algebra over a
field $\k$ of characteristic zero. In this section we see the
following: if $A$ is generated by $A_1$, the subspace of degree
$1$ elements, and satisfies few good properties, then the algebra
$A$ can be described as a quotient of the algebra of
differential operators with constant coefficients on the vector space $A_1$.

We will use this theorem to give a description of the cohomology
subalgebra of a projective spherical variety generated by
elements of degree 2, provided that this algebra satisfies Poincare duality.
This theorem is implicit in the work of
Khovanskii-Pukhlikov \cite{Kh-P}. There is also a version of it
for zero-dimensional local Gorenstein algebras (see \cite[Ex. 21.7]{Eisenbud}).
%\footnote{The conclusion of Theorem
%\ref{thm-comm-alg} is stronger than the conclusion of
%\cite[Exercise 21.7]{Eisenbud}, but the latter is true for the more general class of
%zero-dimensional local Gorenstein algebras.}.

\begin{Th} \label{thm-comm-alg}
Let $A$ be a commutative finite dimensional graded algebra over 
$\k$. Write $A = \bigoplus_{i=0}^{n}
A_i$, where $A_i$ is the $i$-th graded piece of $A$.

Suppose the following conditions hold:
\begin{itemize}
\item[(i)] $A$ is generated, as an algebra, by $A_1$.
\item[(ii)] $A_0 \cong A_n \cong \k$.
\item[(iii)] The bilinear map $A_i \times A_{n-i} \to A_n \cong \k$
given by $(u,v) \mapsto uv$ is non-degenerate for all
$i=0,\ldots,n$.
\end{itemize}
Let $\{ {\bf a_1},\ldots,{\bf a_r} \}$ be a basis for $A_1$ and $P:
\k^r \to \k$ be the polynomial defined by $P(x_1,\ldots,x_n) =
(x_1{\bf a_1}+\cdots+x_r{\bf a_r})^n \in A_n \cong \k$. Then $A$, as
a graded algebra, is isomorphic to the algebra $\k[t_1, \ldots, t_r]
/ I$, where $I$ is the ideal defined by
$$I = \{ f(t_1, \ldots, t_r) \mid
f(\frac{\partial}{\partial{x_1}},\ldots,
\frac{\partial}{\partial{x_r}})\.P = 0 \},$$ and
$f({\partial}/{\partial{x_1}},\ldots, {\partial}/{\partial{x_r}})$
is the differential operator obtained by replacing $t_i$ with
$\partial / \partial{x_i}$ in $f$.
\end{Th}
\begin{proof}[Sketch of proof]
It is easy to see that $I$ is indeed an ideal of
$\k[t_1,\ldots,t_r]$. Consider the surjective homomorphism
$\Phi: \k[t_1,\ldots,t_r] \to A$ given by $f \mapsto f({\bf
a_1},\ldots,{\bf a_r})$. 
%where $f(t_1,\ldots,t_r) =
%\sum_{\beta_1,\ldots,\beta_r}
%c_{\beta_1,\ldots,\beta_r}{t_1}^{\beta_1}\cdots{t_r}^{\beta_r}$ is
%a polynomial in $t_1,\ldots,t_r$ and $ f({\bf a_1},\ldots,{\bf
%a_r}) = \sum_{\beta_1,\ldots,\beta_r} c_{\beta_1,\ldots,\beta_r}{\bf
%a_1}^{\beta_1}\cdots{{\bf a_r}}^{\beta_r}.$ 
We wish to prove that $\ker(\Phi) = I$ and hence $A \cong \k[t_1,\ldots,t_r] / I$.
For this one notes that from the definition, both $I$ and
$\ker(\Phi)$ are generated by homogeneous elements. Take a homogeneous 
polynomial $f$ and consider two cases: 
a) $\deg(f) = n$, b) $\deg(f) < n$. In each case, it is straight forward to verify
that $f \in I$ if and only if $f \in \ker(\Phi)$. In the 
latter case we need to use the assumption (iii) (Poincare duality).
\end{proof}

In Section \ref{sec-main}, we take $A$ to be the cohomology subalgebra generated by
{the Chern classes of line bundles} where $X$ is a projective spherical variety.
For each $i$, the $i$-th graded piece $A_i$ will be $A
\cap H^{2i}(X, \r)$.

\section{Degree of a spherical variety} \label{sec-degree}
Let $G$ be a connected reductive algebraic group over $\c$ and $X$ a
projective spherical $G$-variety of dimension $n$. Recall
that a normal $G$-variety is spherical if a Borel subgroup of $G$ has a
dense orbit. Equivalently, $X$ is spherical if and only if for any
$G$-linearized line bundle $L$ on $X$, the decomposition of the
space of global sections $H^0(X, L)$ into irreducible $G$-modules is multiplicity free.
In this section, we explain the formulae for the degree
of a $G$-linearized line bundle $L$ on $X$ in terms of: 1) an integral over the
moment polytope $\mu(X, L)$, and 2) the volume of a larger polytope $\Delta_{\s}(X, L)$.
We will state the theorems for very ample
line bundles although since $\deg(L^{\otimes k}) = k^n \deg(L)$,
it is easily seen that the formulae hold for ample line bundles as well. \\

\noindent{\bf Notation:} Throughout the rest of the paper we will
use the following notation. $G$ denotes a connected reductive algebraic group (over $\c$),
$B$ a Borel subgroup of $G$, and $T$ a maximal torus contained in $B$. The Lie algebra of $T$ and
its dual will be denoted respectively by $\t$ and $\t^*$. Let $\Phi
= \Phi(X, T) \subset \t^*$ be the root system with $\Phi^+ =
\Phi^+(X, T)$ the subset of positive roots. We denote by $\alpha_1,
\ldots, \alpha_r$ the corresponding simple roots where $r$
is the semi-simple rank of $G$. Let $W$ be the Weyl group of $(G,T)$,
and $s_1, \ldots, s_r \in W$ the simple reflections associated with $\alpha_1,
\ldots, \alpha_r$. These generators of $W$ define a length function
$\ell$ on this group. We denote by $w_0$ the unique longest element in
$W$. Then $N = \ell(w_0)$ is the number of positive roots. Denote by
$\Lambda$ the weight lattice of $G$ (that is, the character group of
$T$), and by $\Lambda^+$ the semigroup of dominant weights. Put
$\Lambda_\r = \Lambda \otimes_{\z} \r$. Then the convex cone
generated by $\Lambda^+$ in $\Lambda_\r$ is the positive Weyl
chamber $\Lambda^+_{\r}$. For a dominant weight $\lambda \in
\Lambda$, the irreducible $G$-module with highest weight $\lambda$
will be denoted by $V_\lambda$, and $v_\lambda$ will represent a highest
weight vector for $V_\lambda$.

\subsection{Formula for degree} \label{subsec-deg-sph}
Let $L$ be a very ample $G$-linearized line bundle on the projective spherical
variety $X$. It gives rise to an embedding $X \hookrightarrow
\p(V)$, with $V = H^0(X, L)^*$. In this section we state
the Brion-Kazarnovskii formula for the degree of $X$ in $\p(V)$.

Recall that for every projective $G$-variety $X$ with a $G$-linearized line
bundle $L$, the moment polytope is defined as:
\begin{equation} \label{moment-conv-hull}
\mu(X, L) = \overline{\textup{conv}(\bigcup_{k}\{\lambda/k \mid
\textup{multiplicity of } V_\lambda^* \textup{ in } H^0(X,
L^{\otimes k}) \neq 0\})}.
\end{equation}

\begin{Rem} \label{rem-moment-symp}
Interestingly, when $X$ is smooth, $\mu(X, L)$ coincides with
the moment polytope (also called the {\it Kirwan polytope}) in symplectic geometry \cite{Brion1}.
%since $L$ is very ample we have a $G$-equivariant embedding
%$X \hookrightarrow \p(V)$, $V = H^0(X, L)^*$. Let $K$ be a maximal compact subgroup
%of $G$ and choose a $K$-invariant hermitian inner product
%on $V$. It induces a Kaehler metric on $\p(V)$ and hence on $X$. The imaginary part
%of this metric gives a $K$-invariant symplectic structure on $X$. With this, $X$
%becomes a $K$-Hamiltonian manifold. Let $\phi_L: X \to \Lie(K)^*$ be its moment map.
%Then $$\mu(X, L) = \phi_L(X) \cap \Lambda^+_\r.$$
\end{Rem}

\begin{Th}[$\S$4.1, \cite{Brion2}] \label{thm-Brion}
The degree of $X$ in $\p(V)$ is equal to
$$ n! \int_{\mu(X, L)} F(\gamma) d\gamma,$$ where $F(\gamma)$ is the
function on the linear span of $\mu(X, L)$ defined by
$$ F(\gamma) = \prod_{\alpha \in \Phi^+ \setminus E}
\frac{\langle \gamma, \alpha \rangle}{\langle \rho,
\alpha\rangle}.$$ Here $\langle \cdot, \cdot\rangle$ is the Killing
form, $E$ is the set of positive roots which are orthogonal to $\mu(X, L)$
and $\rho$ is half the sum of positive roots.
\end{Th}
\subsection{String polytopes} \label{sec-string-polytope}
For the sake of completeness, in this section we recall the
definition of string polytopes for a reductive group $G$. Some of
this recollection is taken from \cite[$\S 1$]{Alexeev-Brion}.

Consider the algebra
$$ A = \c[G]^U$$ of regular functions on G which are invariant under the
right multiplication by $U$. The group $G \times T$ acts on $A$
where $G$ acts on the left and $T$ acts on the right, since it
normalizes $U$. It is well-known that we have the following
isomorphism of $G \times T$ modules $$ A \cong \bigoplus_{\lambda
\in \Lambda^+} V_\lambda^*,$$ where $V_\lambda^* = V_{\lambda^*}$ is
the irreducible representation with highest weight $\lambda^* = -w_0(\lambda)$, and
$T$ acts on each $V_\lambda^*$ via the character $\lambda$.

The vector space $A$ has a remarkable basis $\B = (b_{\lambda,
\phi})$, usually called the {\it dual canonical basis}, such that each $b_{\lambda, \phi}$ is an eigenvector of $T \times T \subset G \times T$, of weight $\lambda$ for the right
$T$-action. For fixed $\lambda$, the vectors $b_{\lambda,\phi}$ form
a basis for $V_\lambda^*$. 
%This basis is the dual basis for the
%so-called {\it canonical basis} of $V_\lambda$
%which consists of the nonzero $bv_\lambda$,
%where $b$ belongs to the specialization of Kashiwara-Luszig's {\it
%canonical basis} at $q=1$.
%The basis $\B$ is called the {\it dual
%canonical basis}. 
For a reduced decomposition of $w_0$, the
longest element of $W$, one defines a parametrization of $\B$ {by integral points},
called a {\it string parametrization} \cite{Littelmann, B-Z}.
Recall that an $N$-tuple of simple reflections
$$\s = (s_{i_1}, s_{i_2}, \ldots, s_{i_N})$$ is a {\it reduced
decomposition} for $w_0$ if $w_0 = s_{i_1}s_{i_2}\cdots s_{i_N},~N =
\ell(w_0)$. The string parametrization associated to $\s$ is an
injective map
$$\iota_{\s}: \B \to \Lambda^+ \times \n^N,$$
$$ \iota(b_{\lambda,\phi}) = (\lambda, t_1, \ldots, t_N).$$
The string parameters have to do with the weight of a basis element
as an eigenvector for the $T\times T$-action: the weight of
$b_{\lambda, \phi} \in \B$ for the left $T$-action is $$-\lambda +
t_1\alpha_{i_1} + \cdots + t_N\alpha_{i_N}.$$

A remarkable property of the string parameterization is that its image consists of all the 
integral points in a certain rational convex polyhedral cone 
$\mathcal{C}$ in $\Lambda_{\r} \times \r^N$ (\cite{Littelmann}).
\begin{Def} The string polytope $\Delta_{\s}(\lambda)$ is the
polytope in $\r^N$ obtained by slicing the cone $\mathcal{C}$ at
$\lambda$, that is $$\Delta_{\s}(\lambda) = \{(t_1, \ldots, t_N)
\mid (\lambda, t_1, \ldots, t_N) \in \mathcal{C} \}.$$
\end{Def}
Note that: 1) $\Delta_{\s}(\lambda)$ is defined for any $\lambda \in \Lambda_\r^+$.
2) From the definition it follows that $\Delta_\s(k\lambda) =
k\Delta_\s(\lambda)$ for any positive integer $k$. 

By what was said
above, when $\lambda$ is a dominant weight,
the lattice points in $\Delta_{\s}(\lambda)$, i.e. the points
in $\Delta_{\s}(\lambda) \cap \z^N $, are in
bijection with the elements of the basis $b_{\lambda, \phi}$ for
$V_\lambda^*$ (and hence in bijection with the basis for
$V_\lambda$). Thus,
\begin{equation} \label{equ-dimV_lambda} \#(\Delta_{\s}(\lambda)
\cap \z^N) = \dim(V_\lambda).
\end{equation}

Let $v_\lambda \in V_\lambda$ be a highest weight vector and
$P_\lambda$ the parabolic subgroup associated to the
weight $\lambda$, that is, $P_\lambda$ is the stabilizer of the
point $[v_\lambda] \in \p(V_\lambda)$. Then we have an embedding
$i: G/P_\lambda \hookrightarrow \p(V_\lambda)$, given by
$gP_\lambda \mapsto g \cdot [v_\lambda]$. Let
$L_\lambda = i^*(\mathcal{O}(1))$ be the line bundle on $G/P_\lambda$ induced by this embedding. 
By Borel-Weil one knows that
for $k>0$, $H^0(X, L_\lambda^{\otimes k}) \cong V_{k\lambda}^*$ as $G$-modules.
Put $n=\dim(G/P_\lambda)$. From
(\ref{equ-dimV_lambda}) it follows that the degree of $G/P_\lambda$, as a subvariety of
$\p(V_\lambda)$ is equal to $n!\Vol_n(\Delta_{\s}(\lambda))$.

\begin{Rem} \label{rem-GC}
{Let $G$ be of type $A$. In the classic paper \cite{G-C}, for each dominant weight $\lambda$ the authors construct 
a polytope $\Delta(\lambda)$ defined by certain explicit inequalities. These are now referred to as the {\it Gelfand-Cetlin polytopes}. 
This construction has been extended to other groups with classical simple Lie algebras in \cite{B-Z}.} 
The Gelfand-Cetlin polytopes are special cases of the string
polytopes. More precisely, let $G=\GL(n, \c)$. The Weyl group
is $W = S_{n+1}$. Let us take the {\it nice} reduced decomposition
$$w_0 = (s_1)(s_2s_1)(s_3s_2s_1) \cdots (s_ns_{n-1}\cdots s_1)$$
for $w_0$, where $s_i$ denotes the transposition exchanging $i$ and
$i+1$. Then $\Delta_{\s}(\lambda)$ can be identified with the
Gelfand-Cetlin polytope corresponding to $\lambda$. Similarly, when
$G = \SP(2n, \c)$ or $\SO(n, \c)$, for a similar choice of a reduced
decomposition, one can recover the Gelfand-Cetlin polytopes as the
string polytopes \cite{Littelmann}.
\end{Rem}

\subsection{Degree as volume of a polytope}
\label{sec-polytope}
The formula for the degree of $(X, L)$ in Theorem \ref{thm-Brion}
can be interpreted as the volume of a certain polytope $\Delta_{\s}(X, L)$.
Fix a reduced decomposition $\s$ for
$w_0$. Define $\Delta_{\s}(X, L)$ to be the
polytope in $\Lambda_{\r} \times \r^N$ such that the first projection
$p:\Lambda_{\r} \times \r^N \to \Lambda_{\r}$ maps $\Delta_{\s}(X,
L)$ onto the moment polytope $\mu(X, L)$ and the fibre of $p:
\Delta_{\s}(X, L) \to \mu(X, L)$ over $\lambda \in \mu(X, L)$ is the
string polytope $\Delta_{\s}(\lambda)$. $\Delta_\s(X, L)$ is a convex rational polytope. We call it
the {\it Newton polytope} of $(X, L, {\s})$.
The Newton polytope has the property that for all non-negative integers $k$:
\begin{equation} \label{equ-Ehrhart}
\#(k\Delta_{\s}(X, L) \cap (\Lambda \times \z^N)) = \dim H^0(X, L^{\otimes k}).
\end{equation}
That is, the {\it Ehrhart function} of $\Delta_{\s}(X,
L)$ coincides with the Hilbert function of $(X, L)$. As in
Section \ref{sec-string-polytope}, it follows that
\begin{Th}[Remark 3.9
(iii)\cite{Alexeev-Brion}] \label{thm-Alexeev-Brion} Let $X$ be a
projective spherical variety of dimension $n$ and $L$ a very ample
$G$-linearized line bundle on $X$. Then
$$\deg(X, L) = n! \Vol_n(\Delta_{\s}(X, L)),$$
where $\Vol_n$ is the $n$-dimensional Lebesgue measure in the subspace spanned by
the Newton polytope $\Delta_\s(X, L) \subset \Lambda_\r \times \r^N$ and normalized with
respect to the lattice $\Lambda(X) \times \z^N$, i.e. the smallest nonzero volume of a parallelepiped with vertices
in the lattice $\Lambda(X) \times \z^N$ is $1$. Here $\Lambda(X)$ is the weight lattice of the spherical
variety $X$.
%(If $\dim(\Delta_\s(X, L)) < n$, its $n$-dimensional volume is $0$.)
\end{Th}

\section{Additivity of moment and string polytopes} \label{sec-additive}
Let $X$ be a projective spherical variety, and $L$ a very ample $G$-linearized line bundle on it.
In this section we investigate the
additivity of maps which associate to $L$ its moment polytope $$L \mapsto \mu(X, L),$$
and its Newton polytope
$$L \mapsto \Delta_{\s}(X, L).$$
(As before $\s$ is a fixed reduced decomposition for the longest element $w_0$.)
Namely, if $L_1$ and $L_2$ are two very ample $G$-linearized line bundles on $X$ we wish to investigate
whether
\begin{equation} \label{equ-additive-moment-polytope}
\mu(X, L_1 \otimes L_2) = \mu(X, L_1) + \mu(X, L_2),
\end{equation}
or
\begin{equation}
\Delta_{\s}(X, L_1 \otimes L_2) = \Delta_{\s}(X, L_1) + \Delta_{\s}(X, L_2),
\end{equation}
where the addition in the righthand side is the Minkowski sum of convex polytopes.
We will see that the additivity of the moment polytope map, and consequently
the Newton polytope map, in general is false. But we will show that it holds
in some important cases: 1) toric varieties, and 2) complete symmetric varieties
and line bundles which restrict trivially to the open orbit. (This in particular includes the group compactifications.)
The additivity of the moment polytope for
toric varieties has been known, while the corresponding result for
symmetric varieties seems to be new.

For a $G$-variety $X$, we denote the group of all $G$-linearized line bundles on $X$ by $\Pic_G(X)$.

\begin{Th}[Additivity of moment polytope for toric varieties] \label{th-toric-additive}
When $X$ is a projective toric variety the moment polytope map is additive
(in the sense of (\ref{equ-additive-moment-polytope})).
\end{Th}
\begin{proof}
Let $\Sigma$ be the fan of the toric variety $X$ with $\Sigma(1)$ the collection of
one-dimensional cones in $\Sigma$. For $\rho \in \Sigma(1)$ let $u_\rho$ denote the smallest
integral vector along $\rho$ (i.e. a primitive vector) and pointing in the opposite
direction. Also let $D_\rho$ denote the irreducible divisor which is the closure of the codimension $1$ orbit
corresponding to $\rho$. It is well-known that every (Cartier) divisor $D$ on $X$ is equivalent to
a linear combination $\sum_{\rho \in \Sigma(1)} a_\rho D_\rho$. The divisor
$D = \sum_{\rho \in \Sigma(1)} a_\rho D_\rho$ is
$T$-invariant and its corresponding line bundle $L(D)$ has a $T$-linearization. Moreover,
if $D$ is very ample then the moment polytope of the $T$-linearized
line bundle $L = L(D)$ is given by the following inequalities:
\begin{equation}
\mu(X, L) = \{ m \mid \langle m , u_\rho \rangle \geq -a_\rho, \forall \rho \in \Sigma(1) \}.
\end{equation}
Now for $i=1,2$, let $D_i = \sum_{\rho \in \Sigma(1)} a_{\rho,i} D_\rho$ be two very ample divisors with the corresponding
line bundles $L_i =L(D_i)$ and moment polytopes
$\mu_i = \mu(X, L(D_i)) = \{ m \mid \langle m , u_{\rho} \rangle \geq -a_{\rho,i}, \forall \rho \in \Sigma(1)\}$.
Then $D_1 + D_2 = \sum_{\rho \in \Sigma(1)} (a_{\rho,1} + a_{\rho,2})D_\rho$ and
$L = L(D_1 + D_2) = L_1 \otimes L_2$. Also the moment polytope $\mu(X, L)$ is given by
$\{ m \mid \langle m , u_\rho \rangle \geq -(a_{\rho,1} + a_{\rho,2}), \forall \rho \in \Sigma(1) \}$
which is clearly the Minkowski sum of the polytopes $\mu_1$ and $\mu_2$. This proves the theorem.
\end{proof}

The degree map $$L \mapsto \deg(X, L) = c_1(L)^n,$$ extends to a polynomial of degree $n$ on the
vector space $\Pic(X)\otimes_\z \r$ (and hence on $\Pic_G(X) \otimes_\z \r$).
%A basic example in which the moment polytope is additive is the class of toric varieties.
The additivity in the toric case leads to a nice description of
the degree function on $\Pic_T(X) \otimes \r$ as follows:
The semigroup of convex polytopes in $\r^n$ with respect to the Minkowski sum has
the cancelation property. Hence it can be extended to a (real) vector space
consisting of {\it virtual polytopes} which are the
formal differences of convex polytopes. Let $\mathcal{V}$ denote this vector space. 
The volume as a function on the space of polytopes
extends to a polynomial function on the vector space $\mathcal{V}$.
Thus the volume function is also called the {\it volume polynomial}.
Let us say that two polytopes $\Delta_1$, $\Delta_2$ are {\it equivalent}, and
write $\Delta_1 \sim \Delta_2$, if $\Delta_1$ can be identified with $\Delta_2$ by a translation
in $\r^n$. It is easy to verify that the equivalence relation $\sim$ extends to an equivalence relation on the vector space
$\mathcal{V}$ respecting the vector space operations.

Let $X$ be a $T$-toric variety.
The moment polytope $\mu(X, L)$ depends on the $T$-linearization of $L$. For
different $T$-linearizations of $L$ the moment polytopes are equivalent.
By additivity, the map $L \mapsto \mu(X, L)$ induces a linear map
$\mu: \Pic(X) \otimes \r \to \mathcal{V}/\sim$. Now from the Kushnirenko theorem,
$$\deg(X, L) = n! \Vol_n(\mu(X,L)).$$
It follows that the degree polynomial $$\deg: \Pic(X) \to \z, \quad L \mapsto \deg(X, L) = c_1(L)^n$$
is the composition of the linear map $\mu$ with the volume polynomial.

The following is an example of a spherical variety for which
the additivity of the moment polytope fails.\footnote{This example
as well as the additivity of the moment polytope for symmetric varieties was suggested to the
author by M. Brion.}
\begin{Ex} \label{ex-Brion}
Let $G = \SL(2,\c)$ act on $X = \c\p^1 \times \c\p^1$ diagonally.
This is a spherical variety. Let $V_n$ denote the irreducible
representation of $G$ of dimension $n+1$ and let $v_n$ be a highest
weight vector of $V_n$. The flag variety $G/B \cong \c\p^1$ embeds
in $\p(V_n)$ via $gB \mapsto g\cdot v_n$. Let $L_n$ be the line
bundle on $\c\p^1$ induced by this embedding. Then $L_n =
\mathcal{O}(n)$. Let $\pi_1$ and $\pi_2$ denote the projection of
$X$ on the first and second factors respectively. For non-negative
integers $a$ and $b$ define
$$L_{a,b} = \pi_1^*(L_a) \otimes \pi_2^*(L_b).$$
It is a $G$-linearized line bundle on $X$ for the diagonal action of $G$. Let
us compute the moment polytopes $\mu(X, L_{2,1})$, $\mu(X, L_{1,2})$
and $\mu(X, L_{3,3})$ to show that the latter is not the sum of the
other two moment polytopes. One has
\begin{eqnarray*}
H^0(X, (L_{2,1})^{\otimes k}) &=& H^0(\c\p^1 \times \c\p^1, \pi_1^*(L_{2k})
\otimes \pi_2^*(L_k)) \cr &=&  H^0(\c\p^1, \pi_1^*(L_{2k})) \otimes
H^0(\c\p^1, \pi_2^*(L_k)) \cr &=& V_{2k} \otimes V_{k} \cr &=& V_{k}
\oplus V_{k+2} \oplus \cdots \oplus V_{3k}.
\end{eqnarray*}
From (\ref{moment-conv-hull}) we obtain that $\mu(X, L_{2,1})$ is
the line segment $[1,3]$. Similarly $\mu(X, L_{1,2})$ is $[1,3]$.
Finally
\begin{eqnarray*}
H^0(X, (L_{3,3})^{\otimes k}) &=& H^0(\c\p^1 \times \c\p^1, \pi_1^*(L_{3k})
\otimes \pi_2^*(L_{3k}))\cr &=& H^0(\c\p^1, \pi_1^*(L_{3k})) \otimes
H^0(\c\p^1, \pi_2^*(L_{3k})) \cr &=& V_{3k} \otimes V_{3k} \cr &=&
V_{0} \oplus V_{2} \oplus \ldots \oplus V_{6k},
\end{eqnarray*}
and thus the moment polytope of $(X, L_{3,3})$ is the line segment
$[0, 6]$. It is obvious that, regarded as $1$-dimensional polytopes,
the Minkowski sum of $[1,3]$ with itself is not equal to $[0, 6]$.
\end{Ex}

\begin{Rem}
\noindent 1)
{Let $\textup{Amp}_G(X)$ denote the cone in $\Pic_G(X)$ generated by the 
ample $G$-linearized line bundles on $X$. For a projective spherical variety $X$ this cone
is finitely generated. One can show that there is a subdivision of the cone
$\textup{Amp}_G(X)$  into smaller cones such that the moment polytope map $L \mapsto \mu(X, L)$ is additive restricted to 
each cone of the subdivision. In the above example 
$L_{a,b} \mapsto \mu(X, L_{a,b}) =[|a-b|, a+b]$ is additive on the cones $a \leq b$ and 
$a \geq b$ respectively. }\\

\noindent 2)
{Suppose $\sigma$ is a cone of maximal dimension in $\textup{Amp}_G(X)$ on which $L \mapsto \mu(X, L)$ 
is additive. Extend the map $\mu(X, \cdot)_{|\sigma}$ to the whole $\Pic_G(X)$ by linearity and call it 
$\tilde{\mu}(X, \cdot)$. The moment polytope $\mu(X, L)$ may not coincide with $\tilde{\mu}(X, L)$ for every $L \in \textup{Amp}_G(X)$, but still the degree formula in Theorem \ref{thm-Brion} holds for $\tilde{\mu}(X, L)$ instead of $\mu(X, L)$. }
\end{Rem}

Now we prove the additivity of the moment polytope
map for a complete symmetric variety $X$ and line bundles which are trivial
restricted to the open orbit. One can verify that in the above example the
line bundles $L_1$, $L_2$ restrict to non-trivial line bundles on the open orbit. In other words,
$L_1$ and $L_2$ do not correspond to any $G$-invariant divisor.
For simplicity we now assume $G$ is semi-simple adjoint.
Although most of the discussion in the rest of this section, with slight
modification, holds for any connected reductive group.

A homogeneous space $G/H$ is a {\it symmetric homogeneous space} if $H$ is
the fixed point set of an algebraic involution $\sigma$ of G (i.e. an
algebraic automorphism of order $2$). A normal $G$-variety which
contains $G/H$ as an open orbit is called a {\it symmetric variety}.
Symmetric varieties are spherical.

A torus in $G$ is $\sigma$-split if for every $x$ in the torus $\sigma(x) = x^{-1}$. All
maximal $\sigma$-split tori are conjugate. It is well-known that one can
choose a maximal $\sigma$-split
torus $T_1$ and a maximal torus $T$ which is $\sigma$-invariant and contains $T_1$.
The action of $\sigma$ on the roots $\Phi = \Phi(G, T)$ gives a decomposition
$\Phi = \Phi_0 \cup \Phi_1$ where $\Phi_0$ are the $\sigma$-fixed roots and
$\Phi_1$ are the roots which are moved under $\sigma$. Fix a set of positive
roots $\Phi^+$ such that $\sigma(\Phi^+ \cap \Phi_1) \subseteq \Phi^- \cap \Phi_1$.
Let $\Delta$ be the set of simple roots corresponding to this set of positive roots.
Again write $\Delta = \Delta_0 \cup \Delta_1$ as fixed and moved roots. The elements
of $\Phi_1$ restricted to $\Lie(T_1)$ constitute a root system with basis
the restriction of the elements in $\Delta_1$ to $\Lie(T_1)$. Let $W_1$ be the
Weyl group of this restricted root system.

The variety $G/H$ has a distinguished point, namely $x_0 = eH$, the coset of identity.
The $T_1$-stabilizer of $x_0$ is $T_1 \cap H$ which is finite.
Put $S = T_1 / (T_1 \cap H)$. Then  $S$ is a torus which acts with no stabilizer at $x_0$.
We can write $\Lie(T) = \Lie(T_1) + \Lie(T_1)^{\bot}$ where $\Lie(T_1)^{\bot}$ is the
orthogonal complement to $\Lie(T_1)$ under the Killing form. The orthogonal projection $\Lie(T)
\to \Lie(T_1)$ gives an embedding $\Lie(T_1)^* \hookrightarrow \Lie(T)^*$ with image
$\textup{ann}(\Lie(T_1)^{\bot})
=\{ \xi \in \Lie(T_1) \mid \langle \xi, t \rangle = 0,~\forall t \in \Lie(T_1)^{\bot} \}$.
We will identify $\Lie(S)^* = \Lie(T_1)^*$ with this subspace.
It can be verified that under this identification, the weights of $S$ go to the weights of $T$, which
identifies $\Lambda(S)$, the weight lattice of $S$ (respectively its real span $\Lambda(S)_\r$)
with a sublattice of $\Lambda$ (respectively a subspace of $\Lambda_\r$).
{One also verifies that under this identification, the
intersection of each Weyl chamber for $T$ with $\Lambda(S)_\r$ is a Weyl chamber for $W_1$.
Moreover, the positive Weyl chamber for $T$ intersected with $\Lambda(S)_\r$ gives the
positive Weyl chamber $\Lambda_\r^+(S)$ for $W_1$ (for the choice $\Delta_1$ of basis roots).}

Let $X$ be a regular compactification of $G/H$ (in the sense of \cite{DP}). We refer to $X$ as a
{\it complete symmetric variety}. Let $Z = Z_X$ be the closure of $S \cdot x_0$ in $X$. It is a smooth $S$-toric variety.
One shows that the fan of $Z$ in $\Lie(S)$ is $W_1$-invariant and hence $Z$ is
determined by a fan in the positive Weyl chamber $\Lambda(S)_{\r}^+$ for $W_1$.

%\begin{Th} [\cite{DP}]
%$X \mapsto Z$ is a one-to-one correspondence between the $G$-equivariant
%completions of $G/H$ and the complete $S$-toric varieties whose fans are
%$W_1$-invariant. Moreover, $X$ is smooth if and only if $Z$ is smooth.
%\end{Th}
Let $L$ be a very ample $G$-linearized line bundle on $X$. Moreover, assume that
the restriction of $L$ to the open orbit is trivial.
The following result due to Bifet-Deconcini-Procesi shows how the moment maps of $(X, L)$ and $(Z, L_{|Z})$
are related.
\begin{Th}[\S{12}, \cite{BDP}] \label{th-BDP}
1) The polytope $\mu(X, L)$ lies in the cone $\Lambda^+_\r(S)$ (under the above identification of $\Lambda_\r(S)$ with a subspace of
$\Lambda_\r$. 2) The $W_1$-orbit of $\mu(X, L)$ coincides with $\mu(Z, L_{|Z})$.
\end{Th}

\begin{Th}[Additivity of moment polytope for complete symmetric varieties]
\label{th-symm-var-additive}
For a complete symmetric variety $X$ and line bundles which restrict trivially to the open orbit, 
the moment polytope map is
additive (in the sense of \ref{equ-additive-moment-polytope}).
\end{Th}

We need the following lemma about Weyl groups of root systems.
Let $\W$ be a Weyl group of a root system in $\r^r$.
Fix a positive Weyl chamber $C$. It is a simplicial cone.
Let $\alpha_1, \ldots, \alpha_r$ be simple roots with simple reflections $s_{\alpha_1}, \ldots,
s_{\alpha_r}$. Let $\omega_1, \ldots, \omega_r$ be the generators of the simplicial cone $C$.

For $\lambda \in C$, let $P_\lambda$ be the convex hull of $\{ w \cdot \lambda \mid w \in \W\}$.
From the definition it follows that the
polytope $P_{\lambda} \cap C$ is defined by the
inequalities:
\begin{eqnarray*}
\langle x, \alpha_i \rangle &\ge& 0, \quad i=1,\ldots,r \cr \langle
x, \omega_i \rangle &\le& \langle \lambda, \omega_i\rangle, \quad i=1,\ldots,r \cr
\end{eqnarray*}
where $\langle \cdot, \cdot \rangle$ is the standard inner product on $\r^r$.

\begin{Lem} \label{lem-P-lambda}
For $\lambda_1, \lambda_2 \in C$ we have
$$(P_{\lambda_1}\cap C) + (P_{\lambda_2}\cap C) = (P_{\lambda_1 + \lambda_2} \cap C).$$
\end{Lem}
\begin{proof}
It is easy to see that the Minkowski sum of $(P_{\lambda_1}\cap C)$ and
$(P_{\lambda_2}\cap C)$ is given by the inequalities
\begin{eqnarray*}
\langle x, \alpha_i \rangle &\ge& 0, \quad i=1,\ldots,r \cr
\langle x, \omega_i \rangle &\le& \langle \lambda_1, \omega_i \rangle +
\langle \lambda_2, \omega_i \rangle, \quad i=1,\ldots,r
\cr
\end{eqnarray*}
which is the same inequalities defining the polytope
$P_{\lambda_1 + \lambda_2} \cap C$.
\end{proof}

\begin{proof}[Proof of Theorem \ref{th-symm-var-additive}]
Let $L_1$ and $L_2$ be very ample $G$-linearized line bundles on $X$.
Put $\Delta_1 = \mu(X, L_1)$, $\Delta_2 = \mu(X, L_2)$ and
$\Delta = \mu(X, L_1 \otimes L_2)$.
We want to show $\Delta = \Delta_1 + \Delta_2.$
Let $W_1 \cdot \Delta_i = \Delta'_{i},~ i=1,2$ and $W_1 \cdot \Delta = \Delta'$.
By Theorem \ref{th-BDP}, $\Delta'_1 = \mu(Z, {L_1}_{|Z}), \Delta'_2 = \mu(Z, {L_2}_{|Z})$
and $\Delta' = \mu(Z, (L_1\otimes L_2)_{|Z})$. From Theorem \ref{th-toric-additive},
$\Delta' = \Delta'_1 + \Delta'_2$.
Thus we need to show that $(\Delta'_1 + \Delta'_2) \cap
\Lambda_\r^+(S) = \Delta_1 + \Delta_2$. That is, if for $\lambda_i \in
\Delta_i$ and $w_1, w_2 \in W_1$ we have $w_1(\lambda_1) + w_2(
\lambda_2) \in \Lambda_\r^+(S)$ then $w_1(\lambda_1) + w_2(\lambda_2)
\in \Delta_1 + \Delta_2$.
Now $w_i(\lambda_i) \in P_{\lambda_i}$
and hence $w_1(\lambda_1) + w_2(\lambda_2) \in P_{\lambda_1}
+ P_{\lambda_2} = P_{\lambda_1+\lambda_2}$. Thus
$w_1(\lambda_1) + w_2(\lambda_2) \in P_{\lambda_1 +
\lambda_2} \cap \Lambda_\r^+(S)$. From Lemma \ref{lem-P-lambda} applied to the restricted
root system $\Phi_1$ and its Weyl group $W_1$,
we conclude that
$$w_1(\lambda_1) + w_2(\lambda_2) \in (P_{\lambda_1} \cap \Lambda_\r
^+(S)) + (P_{\lambda_2} \cap \Lambda_\r^+(S)).$$ But $\lambda_i$ lies in
$\Delta_i$ and hence $P_{\lambda_i} \subset \Delta_i$ which implies that $P_{\lambda_i} \cap \Lambda_{\r}^+(S) \subseteq \Delta_i$. Thus, $w_1(\lambda_1) + w_2(\lambda_2) \in \Delta_1 +
\Delta_2$ as required.
\end{proof}

\begin{Ex} \label{ex-reg-gp-comp}
An interesting example of a symmetric variety is
a group compactification. Below we describe the moment polytope of a
group compactification.
Let $\pi: G \to
\GL(V)$ be a representation of $G$.
Assume that the map $\tilde{\pi}: G \to \p(\End(V))$ is an embedding.
Such a representation is called {\it projectively faithful}.
Let $X_\pi = \overline{\tilde{\pi}(G)}$ be the closure of the image of $G$
in $\p(\End(V))$. 
It enjoys an action of
$G\times G$ which comes from the standard action of $\pi(G) \times
\pi(G)$ on $\End(V)$ by multiplication from left and right. This
action extends the usual left-right action of $G \times G$ on $G$.
Let $P_\pi$ denote the convex hull of the weights of $\pi$. It is a
$W$-invariant polytope in $\Lambda_\r$ and usually called the {\it weight 
polytope of $\pi$}. 
%Recall that the dual fan
%of a polytope (in $\Lambda_\r^*$) is the fan (in $\Lambda_\r$)
%defined by the following property: for any
%$i$, every $i$-dimensional cone of the fan is perpendicular (with respect
%to the natural paring between $\Lambda_\r$ and its dual) to a
%face of the polytope of codimension $i$. One has
%\begin{Th}
%\begin{itemize}\item[(a)]
%The variety $X_\pi$
%consists of a finite number of $(G \times G)$-orbits. These are in
%one-to-one correspondence with the orbits of $W$ acting on the faces
%of the polytope $P_\pi$.
%\item[(b)]Let $\pi$ and $\rho$ be two projectively faithful representations
%of $G$. The normalizations of the corresponding compactifications
%$X_\pi$ and $X_\rho$ are
%isomorphic if and only if the dual fans of the polytopes $P_\pi$
%and $P_\rho$ coincide. If the dual fan of the first polytope is a
%subdivision of the dual fan of the second polytope then there
%exists an equivariant morphism from the normalization of $X_\pi$ to
%that of $X_\rho$ and vice versa.
%\end{itemize}
%\end{Th}
%For Part (a) see \cite[Proposition 8]{Timashev} or \cite[Theorem 2.4.2]{Kapranov}.
%Part (b) follows from Luna-Vust theory for spherical embeddings.
%See \cite[Theorem 4.1]{Knop} for the general theory and \cite[\S 7 and \S 8]{Timashev} for a
%discussion of Luna-Vust theory for group compactifications.
%In fact, Part (a), in the context of algebraic monoids,
%follows from an earlier result of Putcha \cite{Putcha}, and a dual version of
%first part of Part (b) can be found in \cite{Renner}.
({With some extra assumptions 
on $\pi$, $X_\pi$ is a regular compactification of $G$ and in particular smooth,
in general it may not be smooth or even normal.})

The restriction of $\mathcal{O}_{\p(\End(V))}(1)$ to $X_\pi$ gives a very
ample $(G \times G)$-linearized line bundle $L_\pi$.
The polytope $P_\pi$ and $\mu(X_\pi, L_\pi)$ are
related as follows.
Define $\iota: \t^* \to \t^* \oplus \t^*$ by
\begin{equation} \label{equ-iota}
\iota(x) = (x, x^*),
\end{equation}
with $x^* = -w_0(x)$.
Then the image of $P_\pi$ under $\iota$ intersected with
$\Lambda_\r^+ \times \Lambda_\r^+$, the positive Weyl chamber for $G
\times G$, is the moment polytope $\mu(X, L)$.
%Conversely, let $X$ be a $(G \times G)$-equivariant compactification
%of $G$ and $L$ a very ample $(G \times G)$-linearized line bundle
%on $X$. Then $V = H^0(X, L)^*$ is a $(G \times G)$-module and hence a
%$G$-module where $G$ acts diagonally. It gives rise to a
%projectively faithful representation $\pi: G \to \GL(V)$ and
%$X_\pi$ is isomorphic to $X$ as a $(G \times G)$-variety.
The above description of the moment polytope of a group compactification can be
found in \cite{Kazarnovskii}, also the additivity of the weight polytope was noticed by
Kazarnovskii.
\end{Ex}

We end this section with few words about the additivity of the string
polytopes. Let $G$ be a group whose Lie algebra is a direct sum of (complex) classical simple Lie algebras and/or
a commutative Lie algebra. As in \cite{B-Z} for each $\lambda \in \Lambda^+_\r$ one defines a Gelfand-Cetlin
polytope $\Delta(\lambda)$. For a specific choice of a reduced decomposition $\s$, the Gelfand-Cetline polytopes are string polytopes (see Remark \ref{rem-GC}). One can explicitly write down the defining inequalities of the Gelfand-Cetlin polytopes and from these inequalities, 
the additivity of the map $\lambda \mapsto \Delta(\lambda)$ follows, namely for 
$\lambda_1$, $\lambda_2 \in \Lambda^+_\r$ we have:
$$\Delta(\lambda_1 + \lambda_2) = \Delta(\lambda_1)
+ \Delta(\lambda_2).$$
But in general
this is not true for string polytopes. In fact, for a reduced
decomposition $\s$, there is a fan $\Sigma_{\s}$ such that the
following holds: $\lambda_1, \lambda_2 \in \Lambda^+_{\r}$ lie in the same cone
of $\Sigma_{\s}$, if and only if the polytope
$\Delta_{\s}(\lambda_1+\lambda_2)$ is the Minkowski sum of the polytopes
$\Delta_{\s}(\lambda_1)$ and $\Delta_{\s}(\lambda_2)$ (\cite[Lemma 4.2]{Alexeev-Brion}).

From the above paragraph and Theorem \ref{th-symm-var-additive}
it follows that
\begin{Cor} \label{cor-Newton-polytope-additive}
Let $G$ be a connected reductive group whose Lie algebra is a direct sum of (complex) classical simple Lie algebras and/or a commutative Lie algebra, and $\s$ the nice decomposition giving rise to
the Gelfand-Cetlin polytopes. Then the map $L \mapsto \Delta_\s(X, L)$ is additive in the following two cases:
\begin{itemize}
\item[1)] $X$ is a complete symmetric variety and $L$ restricts to a trivial bundle on the open orbit.
\item[2)] $X$ is a (partial) flag variety. 
\end{itemize}
\end{Cor}

\begin{Rem} \label{rem-non-trivial-on-open-orbit}
{The variety in Example \ref{ex-Brion} is a symmetric variety (the automorphism $\sigma$ can be taken to be
conjugation by a semisimple element with distinct eigenvalues). On the other hand, one verifies that the line bundles $L_{1,2}$ and
$L_{2,1}$ do not restrict to trivial line bundles on the open orbit.}
\end{Rem}

{As in the toric case, the formula for the degree (Theorem \ref{thm-Alexeev-Brion}) and the additivity of 
the Newton polytope (Corollary \ref{cor-Newton-polytope-additive})
imply a Bernstein theorem:} 
\begin{Th} \label{th-Bernstein-for-symm-var}
Let $G$ and $\s$ be as above. Let $X$ be an $n$-dimensional complete symmetric variety for $G$ and 
$L_1, \ldots, L_n$ $G$-linearized very ample line bundles on $X$ which restrict trivially to the open orbit. Then 
the intersection number of the Chern classes $c_1(L_i)$ is equal to  
$n!$ times the mixed volume of the corresponding Newton polytopes $\Delta_\s(X, L_i)$.
\end{Th}

Similar statement holds for the intersection numbers of Chern classes of equivariant ample line 
bundles on a (partial) flag variety.

\begin{Rem} \label{rem-enumerative-application}
{The answers to several enumerative geometry problems can be expressed as 
intersection numbers of divisors in the variety of complete
quadrics in $\c P^n$ which is a complete symmetric variety for $G=\PGL(n, \c)$.
Applying Theorem \ref{th-Bernstein-for-symm-var} to this situation, one can represent the 
solutions to these enumerative problems as mixed volumes of Newton polytopes. As an example, 
we expect that this yields an alternative way to compute the answer to 
a famous problem of Schubert: the number of space quadrics tangent to 
$9$ quadrics in general position is $666,841,088$. In \cite[Section 10]{DP} this 
is computed with a different method.}
\end{Rem}

\section{Main theorem} \label{sec-main}
Let $X$ be a projective spherical variety of dimension
$n$. Also assume that the subalgebra $A$ generated by {the Chern classes of line bundles} (over $\r$)
has Poincare duality. We combine the formula for the degree (Theorem
\ref{thm-Alexeev-Brion}) and Theorem \ref{thm-comm-alg} to give a
description of the subalgebra $A$.

%In Theorem \ref{thm-comm-alg}, let $\k=\r$ and $A$ be the subalgebra
%of the cohomology generated by elements of degree 2.
%By Proposition
%\ref{prop-Poincare-duality}, $A$ satisfies the conditions in the
%statement of Theorem \ref{thm-comm-alg}.
%Let $\mathcal{L}$ denote
%the set of all very ample line bundles on $X$. It is a semigroup in
%$\Pic(X)$ and generates $H^2(X,\r)$ as a vector space.
Fix a basis $\{ L_1, \ldots, L_r\}$ for $\Pic(X) \otimes \r$ consisting of
very ample line bundles. 
%Then ${\bf a_i} = c_1(L_i)$, $i=1, \ldots, r$, is a vector space basis for
%$H^2(X, \r)$.

%{\bf We regard $\Pic(X)$ as a subgroup of $H^*(X, \r)$ via $L \mapsto c_1(L)$, the first
%Chern class of $L$. (Recall that when $X$ is 
%smooth it has a paving by (complex) affine spaces and hence $H^2(X, \r) \cong \Pic(X) \otimes \r$.)}

Let $P(x_1, \ldots, x_r)= (x_1 c_1(L_1) + \cdots + x_r c_1(L_r))^n$, be the
homogeneous polynomial of degree $n$ in Theorem \ref{thm-comm-alg}.
It suffices to know the polynomial $P$ for very ample line bundles $L$. From
Theorem \ref{thm-Alexeev-Brion}, for $L = L_1^{x_1} \otimes \cdots
\otimes L_r^{x_r}$ we have
\begin{eqnarray*}
P(x_1, \ldots, x_r) &=& (x_1c_1(L_1)+\cdots+x_rc_1(L_r))^n \cr &=&
\deg(X,L) \cr &=& n! \Vol_n(\Delta_{\s}(X, L)). \cr
\end{eqnarray*}
(As usual, $\s$ is a fixed reduced decomposition for $w_0$.)
From Theorem \ref{thm-comm-alg} we have:
\begin{Th} \label{thm-main}
Let $X$ be a projective spherical variety and assume that the subalgebra $A$ generated by the 
Chern classes of line bundles (over $\r$) has Poincare duality. Let $x_1, \ldots,
x_r$ and $P$ be as above. Then $A$ is
isomorphic to the algebra $\r[t_1, \ldots, t_r] / I$ where $I$ is
the ideal defined by
$$I = \{ f(t_1, \ldots, t_r) \mid
f(\frac{\partial}{\partial{x_1}},\ldots,
\frac{\partial}{\partial{x_r}})\.P = 0 \}.$$ Here
$f(\frac{\partial}{\partial{x_1}},\ldots,
\frac{\partial}{\partial{x_r}})$ is the differential operator
obtained by replacing $t_i$ with $\partial / \partial{x_i}$ in
$f$.
\end{Th}

\section{Examples} \label{sec-examples}
\subsection{Toric varieties} \label{subsec-toric}
Let $G=T=(\c\setminus\{0\})^n$ be a torus.
%Denote by $M$ and $N$ the lattice of characters and one-parameter subgroups
%of $T$ respectively. There is a natural pairing between $M$ and $N$ which identifies
%$N$ with $M^*$. Put $M_\r = M \otimes_\z \r$, $N_\r = N \otimes_\z \r$.
%Then $M_\r \cong N_\r \cong \r^n$.
Let $X$ be a smooth projective $T$-toric variety. $X$ is determined by a fan $\Sigma = \Sigma(X)$.
Let $\Delta$ be a convex lattice polytope of full dimension $n$ and
normal to the fan $\Sigma$. That is, $\Delta$ has integral vertices and
its facets are orthogonal to the one-dimensional cones in $\Sigma$. Every such polytope
has a unique representation
$$\Delta =  \{ m \mid \langle m , u_\rho \rangle \geq -a_\rho, \forall \rho \in \Sigma(1) \},$$
where as before $u_\rho$ are the primitive vectors along the rays $\rho$ and pointing toward the interior of
$\Delta$. The polytope $\Delta$ then determines a divisor
$D_\Delta = \sum_{\rho \in \Sigma(1)} a_\rho D_\rho$ and its corresponding line bundle $L_\Delta$.
Let us denote by $\mathcal{V}_\Sigma$ the vector space over $\r$ spanned by such polytopes $\Delta$.
The elements of $\mathcal{V}_\Sigma$ are the virtual polytopes which are normal to $\Sigma$.
One can identify $\mathcal{V}_\Sigma$ with $\r^r$ where $r$ is the number of
$1$-dimensional cones in $\Sigma$.
Recall that two polytopes $\Delta, \Delta'$ are said to be equivalent, 
$\Delta \sim \Delta'$, if they can be identified by a
translation. This equivalence naturally extends to virtual polytopes.
One knows that the mapping $\Delta \mapsto L_\Delta$ gives rise to a natural isomorphism
$$\mathcal{V}_\Sigma/\sim \,\cong \Pic(X) \otimes_\z \r.$$
%From the Kushnirenko theorem $$\deg(X, L_\Delta) = n! \Vol_n(\Delta),$$
As before let $P$ be the volume polynomial on the vector space of virtual polytopes.
We denote its restriction to the subspace $\mathcal{V}_\Sigma$ again by $P$.
Since the volume is translation invariant, $P$ induces a polynomial $\tilde{P}$
on the quotient vector space $\mathcal{V}_\Sigma / \sim$.

One easily sees that the algebra of differential operators on $\mathcal{V}_\Sigma / \sim$ quotient
by the annihilator ideal of $\tilde{P}$ is isomorphic to the algebra of differential operators
on $\mathcal{V}_\Sigma$ quotient by the annihilator ideal of $P$. Now since the cohomology ring of a
smooth projective toric variety is generated as an algebra by its Picard group,
from Theorem \ref{thm-main} we recover the following description of $H^*(X, \r)$:
\begin{Th}[Theorem $\S$1.4, \cite{Kh-P}] 
\label{thm-Khovanskii-Pukhlikov} 
With notation as above, the
cohomology algebra $H^*(X, \r)$ is isomorphic to the algebra
$\Sym(\mathcal{V}_\Sigma) / I$, \footnote{
$\Sym(V)$ denotes the symmetric algebra of a vector space $V$, i.e. the dual of the algebra of polynomials on $V$.}
where $I$ is the ideal defined by
$$I = \{ f \mid f({\partial}/{\partial{x}})\.P = 0 \},$$ and $P(x)$ is the
homogeneous polynomial of degree $n$ on $\mathcal{V}_\Sigma$ defined by
$$ P(\Delta) = \Vol_n(\Delta),$$
for any polytope $\Delta$.
\end{Th}
One can also recover the usual description of the cohomology ring of a
toric variety by generators and relations, due to
Danilov-Jurkiewicz, from the above theorem (see \cite{Timorin}).

%The linear relations and stokes theorem.

\begin{Rem}
{
\noindent 1) In fact, Theorem \ref{thm-Khovanskii-Pukhlikov} holds over $\z$ as well ($H^*(X, \z)$ is generated by the 
classes of orbit closures, and each orbit closure is a transverse intersection of codimension $1$ orbit closures).\\

\noindent 2) Also one can show that Theorem \ref{thm-Khovanskii-Pukhlikov} holds (over $\r$) when $X$ is a complete simplicial 
toric variety (the conditions in Theorem \ref{thm-comm-alg} are satisfied for the cohomology
ring (over $\r$) of a complete simplicial toric variety). 
}
\end{Rem}

\subsection{Complete flag variety}
\label{subsec-G/B}
Let $G$ be a connected reductive group and let $X = G/B$
be the complete flag variety of $G$. Put $N = \dim(G/B)$.
We now apply Theorem
\ref{thm-main} to obtain a description of the cohomology ring
of $X$ in terms of volume of string
polytopes. It is interesting
to note that in this way we obtain analogous descriptions
for the cohomology rings of the flag variety and toric varieties.

Let $\lambda$ be a dominant weight and 
$L_\lambda$ its corresponding line bundle on $G/B$. 
It is known that the map $\lambda \mapsto L_\lambda$ extends to an
isomorphism between $\Pic(G/B)$ and $\Lambda$, and consequently
an isomorphism between $H^2(G/B, \r)$ and $\Lambda_\r = \Lambda \otimes
\r$. Fix a reduced decomposition $\s$. Recall that $\deg(X, L_\lambda)= N!
\Vol_N(\Delta_{\s}(\lambda))$ (Section
\ref{sec-string-polytope}). It is well-known that the
cohomology ring $H^*(G/B, \r)$ is generated by $H^2(G/B, \r)$.
From Theorem
\ref{thm-main} we have the following:
\begin{Cor} \label{cor-G/B}
The cohomology ring $H^*(G/B, \r)$ is isomorphic to
$\Sym(\Lambda_\r) / I$, where $I$ is the ideal defined by
$$I = \{ f \mid f(\partial / \partial x)\.P = 0 \},$$ and $P$ is the homogeneous
polynomial of degree $n$ on the vector space $\Lambda_\r$ defined by
$$P(\lambda) = \Vol_N(\Delta_{\s}(\lambda))$$ for any $\lambda \in \Lambda^+$.
{Moreover, when $\Lie(G)$ is a direct sum of (complex) classical simple 
Lie algebras and/or a commutative Lie algebra, and $\s$ the nice decomposition giving rise to
the Gelfand-Cetlin polytopes,
the map $\lambda \mapsto \Delta_{\s}(\lambda)$ gives a linear embedding $\Delta_\s$
from $\Lambda_\r$ to the vector space $\mathcal{V}$ of virtual polytopes
in $\r^N$. Hence the polynomial $P$ above is
the composition of the linear map $\Delta_\s$ with the volume polynomial.}
\end{Cor}

\begin{Rem} For a weight $\lambda$, by the Weyl dimension
formula, we have
$$\#(\Delta_{\s}(\lambda) \cap \z^N) = \dim(V_\lambda) = \prod_{\alpha \in \Phi^+}
\langle \lambda+\rho, \alpha \rangle / \langle \rho, \alpha
\rangle,$$ and hence
$$\Vol_N(\Delta_{\s}(\lambda)) = \lim_{k\to\infty}
\frac{\#(k\Delta_{\s}(\lambda) \cap \z^N)}{k^N} = \prod_{\alpha \in
\Phi^+} \langle \lambda, \alpha \rangle / \langle \rho, \alpha
\rangle,$$ where $\rho$ is half the sum of positive roots.
Notice that for each $\alpha$, the equation $$\langle x, \alpha
\rangle / \langle \rho, \alpha \rangle = 0$$ defines a hyperplane
perpendicular to the root $\alpha$. If $P$ is a product of the
equations of a collection of hyperplanes perpendicular to the roots,
by a theorem of Kostant \cite{Kostant}, the ideal of differential
operators on $\Lambda_\r$ which annihilate $P$, i.e.
$$\{ f \mid f(\partial / \partial x) \cdot P = 0 \}$$ is generated
by the $W$-invariant operators of positive degree. This shows that
our description of the cohomology of $G/B$ as a quotient of the
algebra of differential operators, coincides with Borel's
description namely $$H^*(G/B, \r) \cong {\Sym(\Lambda_\r)} / I,$$
where $I$ is the ideal generated by the homogeneous $W$-invariant polynomials of
positive degree.
\end{Rem}

\begin{Rem}
{
When $G$ is a simply connected group of type $A$ or $C$, 
e.g. $\SL(n, \c)$ or $\SP(2n, \c)$, Borel's description holds for $H^*(G/B, \z)$ (see \cite{Manivel}).
}
\end{Rem}

\subsection{Variety of complete conics} \label{ex-complete-conics}
{Consider the set $\Q$ of all smooth conics in $\c P^2$. The group $G  = \PGL(3, \c)$ acts transitively on $\Q$. The stabilizer of the conic $x^2 + y^2 + z^2 = 0$ (in the homogeneous coordinates) is $H = \PO(3, \c)$ and hence
$\Q$ can be identified with the homogeneous space $\PGL(3, \c)/\PO(3, \c)$. The subgroup $\PO(3, \c)$ is the fixed point set of the involution $g \mapsto (g^t)^{-1}$ of $G$ and hence $\Q$ is a 
symmetric homogeneous space. In particular, $\Q$ is spherical. Let $V$ be the vector space of all quadratic forms in $3$ variables and $V^*$ its dual. The map which assigns to a conic $C$
its homogeneous equation (respectively equation of the dual conic $C^*$ ) gives an embedding of $\Q$ in $\p(V)$ 
(respectively $\p(V^*)$). Let $X$ be the closure of the set of all conics $(C, C^*)$ in $\p(V) \times \p(V^*)$.
It is called the {\it variety of complete conics}. It is well-known that $X$ is a smooth variety of dimension $5$
(see \cite[Theorem 3.1]{DP}). This variety plays an important role in classical enumerative geometry. 

A Borel subgroup of $G$ is the isotropy group of a flag $(x, \ell)$. There are two $B$-orbits in $\Q$ 
of codimension $1$ which can be described as follows: 1) the set of conics containing $x$, 
2) the set of conics tangent to $\ell$. The closures of these $B$-orbits are classically denoted 
by $\mu$ and $\nu$. It is well-known that the cohomology ring $H^*(X, \r)$ is generated 
by the classes of $\mu$ and $\nu$. 

For $m, n \in \n$ there is an embedding of $X$ in 
$\p(V^{\otimes m} \otimes (V^*)^{\otimes n})$. It gives a $G$-linearized line bundle $L_{m,n}$ 
on $X$ with the first Chern class $m\mu + n\nu$. Let $\omega_1$, $\omega_2$ denote the 
fundamental weights of $G$. One computes that the the moment polytope 
$\mu(X, L_{m,n})$ is the quadrelateral with vertices $0$, $(2m+n)\omega_1$, 
$(m+2n)\omega_2$ and $2m\omega_1 + 2n\omega_2$ (see \cite[Section 2.7]{Brion2}).
From Theorem \ref{thm-Brion} the degree of $L_{m,n}$ is equal to:
$$P(m, n) = \deg(L_{m,n}) = 5! \int_{\mu(X, L_{m,n})} F(\lambda) d\lambda,$$
where $F(\lambda) =\frac{ \langle\lambda, \alpha\rangle
\langle\lambda, \beta\rangle\langle\lambda, \gamma\rangle}{ \langle\rho, \alpha\rangle
\langle\rho , \beta\rangle\langle\rho , \gamma\rangle}$ with 
$\alpha, \beta, \gamma$ the positive roots and $\rho = (\alpha+\beta+\gamma)/2$. Alternatively, let 
$\s$ be the reduced decomposition for the longest element in $G$ as in Remark \ref{rem-GC} (which gives rise to 
G-C polytopes). Then 
$$P(m,n) = 5! \Vol(\Delta_\s(X, L_{m,n})),$$ where 
$\Delta_\s(X, L_{m,n})$ is the Newton polytope constructed in Section \ref{sec-polytope}. It is a polytope 
fibered over $\mu(X, L_{m,n})$ with G-C polytopes as fibers. 
From the above description of $\mu(X, L_{m,n})$ and the additivity of the G-C polytopes, it follows easily that both maps:
$$(m, n) \mapsto \mu(X, L_{m,n}),$$
$$(m,n) \mapsto \Delta_\s(X, L_{m,n})$$ are additive (in the sense of Section \ref{sec-additive}).
One computes $P(m,n)$ to be: 
$$P(m,n) = m^5 + 10m^4n + 40m^3n^2 + 40m^2n^3 + 10mn^4 + n^5.$$

Finally Theorem \ref{thm-main} gives the following description of the cohomology ring 
$H^*(X, \r)$:

\begin{Cor} \label{cor-complete-conics}
The cohomology ring $H^*(X, \r)$ is isomorphic to
$\r[m,n] / I$, where $I$ is the ideal defined by
$$I = \{ f \mid f(\partial / \partial m, \partial / \partial n)\. P = 0 \}.$$
\end{Cor}

Corollary \ref{cor-complete-conics} agrees with the previously known description of the cohomology ring of $X$:
By direct computation we can find that there are no 
homogeneous differential operators of degrees $1$ and $2$ 
which annihilate the function $P(m,n)$. Also there are two homogeneous differential operators 
$g(\partial / \partial m, \partial / \partial n)$, $h(\partial / \partial m, \partial / \partial n)$
of degrees $3$ and $4$ respectively
$$g(m, n) = 2m^3 - 3m^2n + 3mn^2 - 2n^3,$$
$$h(m,n) = 4m^4 - 3m^2n^2 + 4n^4,$$
which annihilate $P$. Thus the ideal $I$ is generated by the two polynomials $g$ and $h$.
See \cite[Section 3.7]{Brion3} and \cite{Kleiman}.
}

\vspace{.5cm}
\noindent Kiumars Kaveh\\Department of Mathematics, University of Pittsburgh\\Pittsburgh, PA, USA\\
{\it Email address:} {\sf kaveh@pitt.edu}
\end{document}